\documentclass[11pt,twoside]{article}
\usepackage{amsmath, amsthm, amscd, amsfonts, amssymb, graphicx, color}
\usepackage[bookmarksnumbered, colorlinks]{hyperref} \usepackage{float}
\usepackage{lipsum}
\usepackage{afterpage}
\usepackage[labelfont=bf]{caption}
\usepackage[nottoc,notlof,notlot]{tocbibind} 

\usepackage{lipsum}
\usepackage{fancyhdr}
\theoremstyle{definition}

\makeatletter
\def\th@newremark{\th@remark\thm@headfont{\bfseries}}
\makeatletter

\begin{document}
‎
\begin{center}

{\Large \bf 
RATE OF CONVERGENCE OF WAVELET SERIES BY CESARO MEANS} 

{\bf  Neyaz A Sheikh $^*$\let\thefootnote\relax\footnote{$^*$Corresponding Author}}\vspace*{-2mm}\\
\vspace{2mm} {\small  National Institute of Technology, Srinagar} \vspace{2mm}

{\bf M Ahsan Ali}\vspace*{-2mm}\\
\vspace{2mm} {\small   National Institute of Technology, Srinagar} \vspace{2mm}

\end{center}

\vspace{4mm}

\begin{abstract}
Wavelet frames have become a useful tool in time freqency analysis and image processing. Many authors worked in the field of wavelet frames and obtained various necessary and sufficient conditions. Ron and Shen [17] gave a charactarization of wavelet frames. Benedetto and Treiber [3], Ron and Shen [17] presented different presentations to the wavelet frames. Any function $f \in L^2(R)$ can be expanded as an orthonormal wavelet series and pointwise convergence and uniform convergence of series have been discussed extensively by various authors [9, 17].
In this paper we investigate the pointwise convergence of orthogonal wavelet series in Pringscheim's sense. Furthermore, we investigate Cesaro $|C,1,1|$ summability and the strong Cesaro $|C,1,1|$ summability of wavelet series.
\end{abstract}

{\bf{1. Introduction}}

Wavelets with local support in the time and frequency domains were defined by A. Grossman and J. Morlet	[6] in 1984 in order to analyze seismic data. The prototypes of wavelets, however, can be found in the work of A. Haar [8]. To identify the underlying structure and to generate interesting examples of orthonormal bases for $L_2(R)$. S. Mallat [12] and Y. Meyer developed the framework of multiresolution analysis. P.G. Lemarie and Y. Meyer [10] constructed wavelets in $S(R^n),$ the space of rapidly decreasing smooth functions. J.O.Stromberg [17] developed spline wavelets while looking for unconditional bases for Hardy spaces. G. Battle [2] and P.G. Lemari [11] developed these bases in the context of wavelet theory. Spline wavelets have exponential decay, but only $C^N$ smoothness (for a finite N depending on the order of the associated splines). I. Daubechies [4]constructed compactly supported wavelets with $C^N$ smoothness. The support of these wavelets increased with the smoothness; in general, to have $C^{\infty}$ smoothness, wavelets must have infinite support.
Y. Meyer [14] was among the first to study convergence results for wavelet expansions; he was followed by G. Walter $([25], [26])$. Kon and Repheat [9] gave certain conditions for the convergence of wavelet series. Meyer proved that under some regularity assumptions on the wavelets, wavelet expansions of continuous functions converge everywhere. In contrast to these results, the pointwise convergence results presented here give almost everywhere convergence (and convergence on the Lebesgue set ) for expansions of general $L^p (1 \le p \le \infty)$ functions. More results about convergence of wavelet series on different spaces have been obtained by Firdous A. Shah, Neyaz A. Sheikh [23]. In this paper we investigate the pointwise convergence of orthogonal wavelet series in Pringscheim's sense. Furthermore, we study Cesaro $|C,1,1|$ summability and the strong Cesaro $|C,1,1|$ summability of wavelet series.

Let ${\psi}_j(x):j=1, 2...$ be orthonormal system of real-valued functions on X. We consider the orthogonal series\\
\[\sum\limits_{j=1}^{\infty}c_{j}\psi_{j}(x),\eqno(1.1)\]
where ${c_{j}}$ is sequence of real numbers (so-called co-efficients)satisfying the condition
\[\sum\limits_{j=1}^{\infty}c_{j}^{2}\le \infty.\eqno(1.2)\]

By Riesz-Fischer theorem, there exists a function $f(x)\in L^{2} = L^2(X,\mathcal{F},\mu)$ such that (1.1) is the generalized Fourier series of f(x) with respect to the system ${\psi_{j}(x)}$.
The partial sums
\[
s_{m}(x)=\sum\limits_{j=1}^{m}c_{j}\psi_{j}(x), m=1,2,...
\]
converges to f(x) in $L^{2}$ norm:
$$\lim\limits_{m\to\infty}\int|s_{m}(x)-f(x)|^{2}d\mu(x)=0\eqno(1.3)$$
Here and in the sequel, the integrals are taken over the entire space X.
It is well known that the condition (1.2) does not ensure the pointwise convergence of the partial sums $s_{m}(x)$ to   f(x) as $m\to\infty.$
 If
$$\sum\limits_{j=1}^{\infty}c_{j}^{2}[log(j+1)]^2<\infty,\eqno(1.4)$$
where the logarithm is to be the base 2, then
$$\lim\limits_{m\to\infty}s_m(x)=f(x)	~a.e.,\eqno(1.5)$$
where f(x) is the sum of (1.1) in the $L^2$- norm (see in (1.3)).
We note that the condition (1.4) is not only sufficient, but also necessary in certain cases for the fulfillment of (1.5), as the following theorem of Tandori.
 If
$$|c_1|\ge|c_2|\ge...\ge|c_j|\ge...$$
 and
$$ \sum\limits_{j=1}^{\infty}c_j^2[log(j+1)]^2=\infty.$$
then one can construct an ONS ${\psi_j(x):j=1,2,...}$ on the unit interval [0,1] endowed with the ordinary Lebesgue measure so that the orthogonal series (1.1) diverges at each point $x\in[0,1].$
The Cesaro (C,1) summability of the orthogonal series (1.1) defined by the arithmetic means
$$\sigma_M(x)=\frac{1}{M}\sum\limits_{m=1}^{M}s_m(x), M=1,2,...$$
of the partial sums is gauranteed by a weaker condition than (1.1). The Menshov-Kaczmarz theorem (see, e.g,. [1, Theorem 2.8.1, p.125]) reads as follows: If
$$\sum\limits_{j=1}^{\infty}c_j^2[loglog(j+3)]^2<\infty,\eqno(1.6)$$.
then
$$\lim\limits_{M\to\infty}\sigma_m(x)=f(x)~a.e.\eqno(1.7)\\$$

The following theorem on the strong Cesaro (C,1) summability of orthogonal series (1.1) was proved by Borgen: Under the condition (1.6), we even have
$$\lim\limits_{M\to \infty}\frac{1}{M}\sum\limits_{m=1}^{M}|s_m(x)-f(x)|^2=0~a.e.$$
We note that in the case of single series of numbers, the notion of strong Cesaro $|C, 1|$ summability is due to G. H. Hardy (where it is called $H_2$ summability). Here our main goal in this paper is to prove an analogous result for wavelet orthogonal series on the strong Cesaro $|C,1,1|$ summability. We also note that the following generalization of Borgen's theorem was proved by Tandori [24]: Let $1 \le v_1<v_2<...v_m<...$ be an arbitrary sequence of natural numbers and set
$$\sigma_M({v}:x)=\frac{1}{M}\sum\limits_{m=1}^{M}s_{v_m}(x), M=1,2,...$$
If condition (1.6)is satisfied, then we have
$$\lim\limits_{M\to\infty}\frac{1}{M}\sum\limits_{m=1}^{M}|s_{v_m}(x)-f(x)^2=0~a.e.$$

{\bf{2. Background: Wavelet series}}

Let ${\psi_{j,k}:j,k=1,2,...}$ be an ONS wavelet  on a finite positive measure space $(X,F,\mu).$ We consider the wavelet series
$$\sum\limits_{j=1}^{\infty}\sum\limits_{k=1}^{\infty}c_{j,k}\psi_{j,k}(x),\eqno(2.1)$$
where ${c_{j,k}:j,k=1,2,...}$ are wavelet coefficients satisfying the condition
$$\sum\limits_{j=1}^{\infty}\sum\limits_{k=1}^{\infty}c_{j,k}^2<\infty.\eqno(2.2)$$
By the Riesz-Fischer theorem, there exists a function g(x) in $L^2= L^2(X,F,\mu)$ such that (2.1) is the wavelet series of g(x) with respect to the system ${\psi_{j,k}(x)}$ and the recangular partial sums
$$s_{m,n}(x)=\sum\limits_{j=1}^{\infty}\sum\limits_{k=1}^{\infty}c_{j,k}\psi_{j,k}(x),~m,~n=1,2,...\eqno(2.3);$$
of the wavelet series (2.1) converge to g(x) in $L^2$-norm:
$$\lim\limits_{m,n\to\infty}\int|s_{m,n}(x)-g(x)|^2d\mu(x)=0.\eqno(2.4)$$
It is clear that condition (2.2) does not ensure the pointwise convergence of the rectangular partial sums $s_{m,n}(x) \to g(x) as ~ m,n\to\infty.$ The extension of Redemacher-Menshov theorem proved by a number of authors reads as follows: If
$$\sum\limits_{j=1}^{\infty}\sum\limits_{k=1}^{\infty}c_{j,k}^2[log(j+1)]^2[log(k+1)]^2<\infty,\eqno(2.5)$$
then
$$\lim\limits_{m,n\to\infty}s_{m,n}(x)=g(x)~a.e.,$$
where g(x) occurs in (2.4). We note that condition (2.5) is also the best possible one that guarantees the a.e. convergence of the wavelet series (2.1).
Indeed, as it was proved in [16] that if the double sequence ${c_{j,k}}$ is such that
$$|c_{j,k}|\ge |c_{j_{1},k_{1}}| for~ all ~ 1\le j\le j_1, ~ 1\le k\le k_1.$$
and
$$\sum\limits_{j=1}^{\infty}\sum\limits_{k=1}^{\infty}c_{j,k}^2[log(j+1)]^2[log(k+1)]^2<\infty,\eqno(2.6)$$
then there exists an ONS ${\psi_{j,k}(x)}$ endowed with the plane Lebesgue measure such that the wavelet series (2.1) diverges a.e.
The a.e. Cesaro (C, 1, 1) summability of the wavelet series (2.1) defined by the arithmetic means
$$\sigma_{M,N}(x)=\frac{1}{MN}\sum\limits_{m=1}^{M}\sum\limits_{n=1}^{N}s_{m,n}(x), M,N=1,2,...\eqno(2.7)$$
of the rectangular partial sums, can be guaranteed under a weaker condition than (2.5).
The extension of the Menshov-Kaczmarz theorem reads as follows:
If
$$\sum\limits_{j=1}^{\infty}\sum\limits_{k=1}^{\infty}c_{j,k}^2[log log(j+3)]^2[log log(k+3)]^2<\infty,\eqno(2.8)$$
then
$$\lim\limits_{M,N\to\infty}\sigma_{M,N}(x)=g(x)~a.e.\eqno(2.9)$$
A double series of real (or complex) numbers
$$\sum\limits_{j=1}^{\infty}\sum\limits_{k=1}^{\infty}u_{j,k}\eqno(2.10)$$
is said to converge in Pringsheim's sense to the sum s if for every $\epsilon > 0$ there exists $K = K(\epsilon)\in \mathcal{N}$ such that
$$|\sum\limits_{j=1}^{m}\sum\limits_{k=1}^{n}u_{j,k}-s|< \epsilon ~if~\min(m,n) > K.$$
Hardy [7]introduced the notion of regular convergence as follows. The double series (2.10) is said to converge regularly if it converges in Pringsheim's sense and, in addition, if each of it's so-called row (single)subseries defined by
$$\sum\limits_{j=1}^{\infty}u_{j,k}$$ converges for each fixed $k\in N,$ as well as each of its so-called column (single)subseries defined by
$$\sum\limits_{k=1}^{\infty}u_{j,k}$$ converges for each fixed $j\in N.$
Without knowing Hardy's definition, the present order rediscovered the notion of regular convergence (where it was called as ``convergence in a restricted sense'') as follows: for every $\mathcal{E}>0$ there exists $\mathcal{K}=\mathcal{K}(\mathcal{E})\in \mathcal{N}$ such that
$$|\sum\limits_{j=1}^{\infty}\sum\limits_{k=1}^{\infty}u_{j,k}|< \epsilon ~ if ~ \max(j_1,k_1)> K and ~ 1\le j_1\le j_2,~ 1\le k_1\le k_2.$$
It is clear that if the double series (2.10) converges absolutely, that is, if
$$\sum\limits_{j=1}^{\infty}\sum\limits_{k=1}^{\infty}|u_{j,k}|<\infty,$$
then it converges regularly, as well: but the converse implication is not true in general.\\

{\bf{3. The Kronecker lemma for wavelet series}}

In the proofs of our Theorems 1 and 2 we will need the extension of the familiar Kronecker lemma from single to double series of numbers. For the reader’s convenience, first we formulate
the Kronecker lemma for single series as follows: If a non-decreasing sequence ${\lambda_j : j = 1, 2,...}$ of positive numbers tends to $\infty$ as $j\to\infty$, and the sequence ${u_j : j = 1, 2,...}$ of real numbers is such that the series
$$\sum\limits_{j=1}^{\infty}\frac{u_j}{\lambda_j}$$ converges, then $$\lim\limits_{m\to\infty}\frac{1}{\lambda_m}\sum\limits_{j=1}^{m}u_j=0.$$\\
We recall that the Kronecker lemma for double series reads as follows:
 Let ${\lambda_{j,k} : j, k = 1, 2,...}$ be
a double sequence of positive numbers satisfying the following conditions:
$$\delta_{1,0}\lambda_{j,k}=\lambda_{j+1,k}-\lambda_{j,k}\ge 0.~\delta_{0.1}\lambda_{j,k}=\lambda_{j,k+1}-\lambda_{j,k}\ge 0.$$
$$\delta_{1,1}\lambda_{j,k}=\lambda_{j+1,k+1}-\lambda_{j+1,k}-\lambda_{j,k+1}+\lambda_{j,k}\ge 0~for all j,k=1,2,...\eqno(3.1)$$
and
$$\lambda_{j,k}\to\infty~as~max(m,n)\to\infty~(or~min(m,n)\to\infty).$$
If the double series
$$\sum\limits_{j=1}^{\infty}\sum\limits_{k=1}^{\infty}\frac{u_{j,k}}{\lambda_{j,k}}~converges ~regularly.\eqno(3.2)$$
where the $u_{j,k}$ are real numbers, then
$$\frac{1}{\lambda_{m,n}}\sum\limits_{j=1}^{m}\sum\limits_{k=1}^{n}u_{j,k}\to 0~as~max(m,n)\to \infty ~(or~min(m,n)
\to \infty).\eqno(3.3)$$\\

{\bf{4.Convergence of wavelet series}}\\

{\bf{Theorem 1.}} {\it If the wavelet series (2.1) is such that condition (2.2) is satisfied, then
$$\frac{1}{MN}\sum\limits_{m=1}^{M}\sum\limits_{n=1}^{N}|s_{m,n}(x) - \sigma_{m,n}(x)|^2 \to 0 ~a.e.~ as ~ max(M,N)\to \infty. \eqno(4.1)$$

Proof. \rm By invoking (2.3)and (2.7), for m, n$\ge2$ we have
$$s_{m,n}(x)-\sigma_{m,n}(x)=\sum\limits_{j=2}^{m}\sum\limits_{k=2}^{n}\frac{(j-1)(k-1)}{mn}c_{j,k}\psi_{j,k}(x),$$
for $m\ge2$ we have
$$s_{m,1}(x)-\sigma_{m,1}(x)=\sum\limits_{j=2}^{m}\frac{j-1}{m}c_{j,1}\psi_{j,1}(x),$$
for $n\ge2$
$$s_{1,n}(x)-\sigma_{1,n}(x)=\sum\limits_{k=2}^{\infty}\frac{k-1}{n}c_{1,k}\psi_{1,k}(x)$$
Due to the orthonormality of the system ${\psi_{j,k}(x)},$ for m, n = 1,2,... we have,
$$\frac{1}{mn}\int|s_{m,n}(x)-\sigma_{m,n}(x)|^2d\mu(x)=\sum\limits_{j=2}^{m}\sum\limits_{k=2}^{n}\frac{(j-1)^2(k-1)^2}{m^3n^3}c_{j,k}^2\eqno(4.2)$$
for $m\ge2$ we have,
$$\frac{1}{m}\int|s_{m,1}(x)-\sigma_{m,1}(x)|^2d\mu(x)=\sum\limits_{j=2}^{\infty}\frac{(j-1)^2}{m^3}c_{j,1}^2;$$
and for $n\ge2$ we have,
$$\frac{1}{n}\int|s_{1,n}(x)-\sigma_{1,n}(x)|^2d\mu(x)=\sum\limits_{k=2}^{\infty}\frac{(k-1)^2}{n^3}c_{1,k}^2$$
Hence it follows that
$$\sum\limits_{m=1}^{\infty}\sum\limits_{n=1}^{\infty}\frac{1}{mn}\int|s_{m,n}(x)-\sigma_{m,n}(x)|^2d\mu(x)$$
~~~~~~~~~~$$=\sum\limits_{j=2}^{\infty}\sum\limits_{k=2}^{\infty}  (j-1)^2(k-1)^2c_{j,k}^2\sum\limits_{m=j}^{\infty}\sum\limits_{j=k}^{\infty}\frac{1}{m^3n^3}$$
~~~~~~~~~~$$+\sum\limits_{j=2}^{\infty}(j-1)^2c_{j,1}^2\sum\limits_{m=j}^{\infty}\frac{1}{m^3}+\sum\limits_{k=2}^{\infty}(k-1)^2c_{1,k}^2\sum\limits_{n=k}^{\infty}\frac{1}{n^3}.\eqno(4.3)$$
where we interchanged the order of summations with respect to j and m, as well as with respect to k and n on the right-hand side in (4.2). For $j\ge 2 ,$ we have,
$$\sum\limits_{m = j}^{\infty}\frac{1}{m^3}\le \int \limits_{j-1}^{\infty}\frac{dt}{t^3}=\frac{1}{2(j-1)^2}.$$
and for $k\ge 2,$ we also have an analogous inequality. Taking into account
these inequalities, the right-hand side in (4.3) does not exceed
$$\frac{1}{4}\sum\limits_{j=2}^{\infty}\sum\limits_{k=2}^{\infty}c_{j,k}^2\le\infty.$$
thanks to (2.2). To sum up, we have proved that
$$\sum\limits_{m=1}^{\infty}\sum\limits_{n=1}^{\infty}\frac{1}{mn}\int|s_{m,n}(x)-\sigma_{m,n}(x)|^2d\mu(x)\le\infty.\eqno(4.4)$$
By the monotone convergence theorem of the Lebesgue integral, it follows from (4.4) that
$$\sum\limits_{m=1}^{\infty}\sum\limits_{n=1}^{\infty}\frac{1}{mn}|s_{m,n}(x)-\sigma_{m,n}(x)|^2\le\infty, a.e.\eqno(4.5)$$
that is, the double series in (4.5) converges absolutely. Consequently, it converges regularly for almost every $x\in  X$. Furthermore, the conditions in(3.1) are clearly satisfied by $\lambda_{j,k} := jk.$ Thus, condition (3.2) is satisfied and we may apply the extension of the Kronecker lemma for the double series in (4.5) (cf. (3.3)) yielding (4.1) to be proved.

Our second new result gives answer to the problem of the strong Cesaro $|C, 1, 1|$ summability of wavelet series.\\

{\bf {Theorem 2.}} {\it  If the double orthogonal series (2.1) is such that condition (2.8) is satisfied, then it is strong Ces`aro $|C, 1, 1|$ summable a.e., that is,
$$\frac{1}{MN}\sum\limits_{m=1}^{M}\sum\limits_{n=1}^{N}|s_{m,n}(x)-g(x)|^2\to 0 ~ a.e. ~ as M,N\to\infty.\eqno(4.6)$$
where g(x) is the sum of (2.1) in the $L^2$-norm (see in (2.4)).}\\

Proof. We start with the elementary inequality
$$|s_{m,n}(x)-g(x)|^2\le2(|s_{m,n}(x)-\sigma_{m,n}(x)|^2+|\sigma_{m,n}(x)-g(x)|^2).$$
hence it follows that for $M,N\ge2$ we have
$$\frac{1}{MN}\sum\limits_{m=1}^{M}\sum\limits_{n=1}^{N}|s_{m,n}(x)-g(x)|^2$$
$$\le\frac{2}{MN}\sum\limits_{m=1}^{M}\sum\limits_{n=1}^{N}|s_{m,n}(x)-\sigma_{m,n}(x)|^2+\frac{2}{MN}\sum\limits_{m=1}^{M}\sum\limits_{n=1}^{N}|\sigma_{m,n}(x)-g(x)|^2.\eqno(4.7)$$
By Theorem 1, the first term on the right-hand side of (4.7) converges to 0 a.e. as $max\{M,N\}\to\infty. $As to the second term there, due to condition(2.8), we may apply the extended Menshov–Kaczmarcz theorem for double orthogonal series to obtain (2.9), hence we conclude that the second term on the right-hand side of (4.7) also converges to 0 a.e. as $M,N\to\infty.$ The proof of (4.6) is complete.\\

We note that the Cesaro $(C, 1, 1)$ summability of the wavelet series $(2.1)$ follows from its strong $Cesaro |C, 1, 1|$ summability. Indeed, by the definition $(2.7)$ and the familiar Cauchy inequality, we may estimate as follows:

\begin{eqnarray*}
|g(x)|~~~\le&&\frac{1}{MN}\sum_{m=1}^{M}\sum_{n=1}^{N}1.|s_{m,n}(x)-g(x)| \\
&&\le\frac{1}{MN}(\sum_{m=1}^{M}\sum_{n=1}^{N}1^2)^{\frac{1}{2}}(\frac{1}{MN}\sum\limits_{m=1}^{M}\sum\limits_{n=1}^{N}|s_{m,n}(x)-g(x)|^2)^{\frac{1}{2}}\\
&&=(\frac{1}{MN}\sum\limits_{m=1}^{M}\sum\limits_{n=1}^{N}|s_{m,n}(x)-g(x)|^2)^{\frac{1}{2}}\end{eqnarray*}
Now, the implication (3.2) $\implies$ (2.9) is clear.

\footnote{This work is supported by DST-SERB Project No.: SR/S4/MS/ 818/13} \\

{\small

\noindent{\bf Neyaz Ahmad Shiekh}

\noindent Department of Mathematics

\noindent Associate Professor of Mathematics

\noindent National Institute of Technology, Srinagar 

\noindent Srinagar, India

\noindent E-mail: neyaz@nitsri.net}\\

{\small
\noindent{\bf  M Ahsan Ali}

\noindent  Department of Mathematics

\noindent Research Scholar

\noindent National Institute of Technology, Srinagar

\noindent Srinagar, India

\noindent E-mail: ahsan$\_$dst$\_$17@nitsri.net}\\

\end{document}